\documentclass{amsart}
\usepackage{amsfonts}
\usepackage{amsmath,amssymb}	
\usepackage{amsthm}
\usepackage{amscd}
\usepackage{graphics}
\usepackage{graphicx}
\theoremstyle{remark}{
\newtheorem{Def}{{\rm Definition}}

\newtheorem{Rem}{{\rm Remark}}
\newtheorem{Prob}{{\rm Problem}}
}

\newtheorem{Prop}{Proposition}
\newtheorem{Thm}{Theorem}

\begin{document}
\title[Reeb graphs of smooth functions on $3$-dimensional manifolds]{On Reeb graphs induced from smooth functions on $3$-dimensional closed manifolds with finitely many singular values}
\author{Naoki Kitazawa}
\keywords{Singularities of differentiable maps; generic maps. Differential topology. Reeb spaces (graphs).\\
\indent {\it \textup{2020} Mathematics Subject Classification}: Primary~57R45. Secondary~57R19.}
\address{Institute of Mathematics for Industry, Kyushu University, 744 Motooka, Nishi-ku Fukuoka 819-0395, Japan\\
 TEL (Office): +81-92-802-4402 \\
 FAX (Office): +81-92-802-4405 \\
}
\email{n-kitazawa@imi.kyushu-u.ac.jp}
\urladdr{https://naokikitazawa.github.io/NaokiKitazawa.html}
\maketitle
\begin{abstract}
The {\it Reeb graph} of a smooth function on a smooth manifold is the graph obtained as the space of all connected components of preimages (level sets) such that the set of all vertices coincides with the set of all the connected components of preimages containing some singular points. 
Reeb graphs are fundamental and important tools in algebraic topological and differential topological theory of Morse functions and their variants.

In the present paper, as a related fundamental and important study, for given graphs, we construct certain smooth functions inducing the graphs as the Reeb graphs. Such works have been demonstrated by Masumoto, Michalak, Saeki, Sharko, among others, and also by the author since 2000s. We construct good smooth functions on suitable $3$-dimensional connected, closed and orientable manifolds.
\end{abstract}

\section{Introduction}
\label{sec:1}
First, throughout the present paper, graphs are fundamental tools. The {\it vertex set} of a graph is the set of all vertices and the {\it edge set} is the set of all edges. A graph is {\it finite} i.e. the vertex set and the edge set are finite and the graph has at least one edge throughout the paper. Such a graph is naturally topologized and a $1$-dimensional polyhedron.

The {\it Reeb graph} (or the {\it Kronrod-Reeb} graph) of a smooth function on a smooth manifold is the graph obtained as the space of all the connected components of preimages (level sets) such that the vertex set coincides with the set of all connected components of preimages containing {\it singular points}: a {\it singular point} of a smooth map is a point at which the rank of the differential of the map is smaller than both the dimensions of the domain and the target. 
Reeb graphs are fundamental and important tools in algebraic topological and differential topological theory of Morse functions and their higher dimensional versions as already seen in \cite{reeb} and so on.

One of fundamental and important studies on Reeb graphs is the following.

\begin{Prob}
\label{prob:1}
For a given graph, can we construct a smooth function (satisfying some conditions) inducing the Reeb graph isomorphic to the graph? 
\end{Prob}

\begin{figure}
\includegraphics[width=35mm]{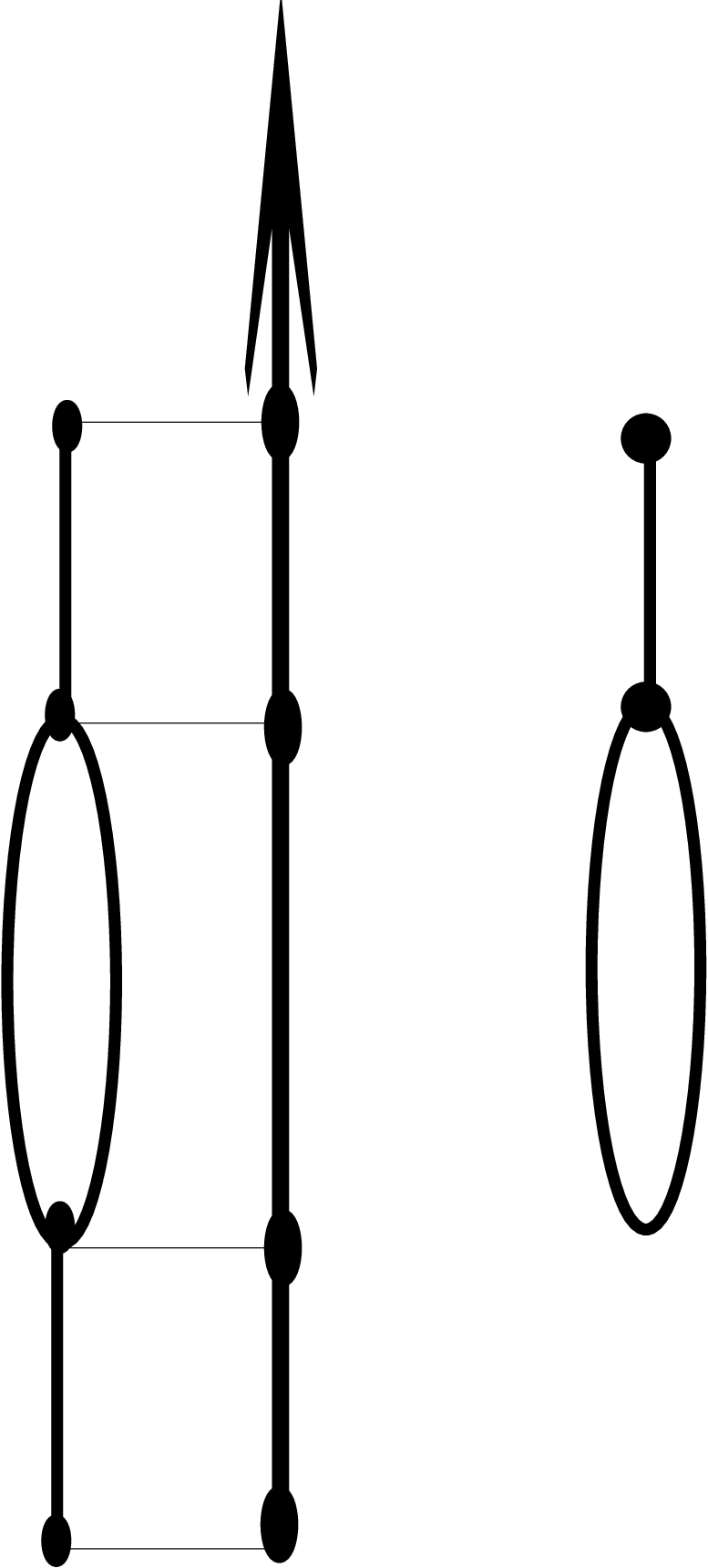}
\caption{A graph admitting a good function and a graph admitting no good function.}
\label{fig:0}
\end{figure}

Construction of a function inducing a Reeb graph isomorphic to a given graph was first demonstrated by Sharko (\cite{sharko}), followed by a work of J. Martinez-Alfaro, I. S. Meza-Sarmiento and R. Oliveira (\cite{martinezalfaromezasarmientooliveira}), one of Masumoto and Saeki (\cite{masumotosaeki}), the works of Michalak (\cite{michalak} and \cite{michalak2}) among others. 
In most of such studies, preimages of {\it regular values}, defined as points in targets we cannot realize as values at some singular points of the functions, are disjoint unions of circles or standard spheres and domains are closed surfaces. The author has first succeeded in explicitly constructing several {\it Morse} functions and so-called {\it circle valued Morse} functions whose preimages of regular values are closed surfaces under several constraints on topologies of preimages of regular values in \cite{kitazawa} and \cite{kitazawa2} (see also Remark \ref{rem:4}). In these works among others and the author, functions inducing Reeb graphs isomorphic to given graphs have been constructed in the following way. First, we prepare a graph admitting a continuous function such that the function obtained by restricting the function to each edge is injective (we call such a function a {\it good function} on the graph): see FIGURE \ref{fig:0}.

It is obvious that a graph has a good function if and only if it has no loops.

Next, we construct local functions around vertices respecting the value of the given good function at each vertex. Last we construct a local function having no singular points around each connected component of the complementary set of the disjoint union of small neighborhoods of vertices.  

In the present paper, we show the following result. Several undefined terminologies and notation will be explained in the next section.
 
\begin{Thm}
\label{thm:1}
Let $G$ be a finite and connected graph $G$ which has at least one edge and has a good function $g:G \rightarrow \mathbb{R}$. Suppose that a non-negative integer is assigned to each edge. Then there exist a $3$-dimensional closed, connected and orientable manifold $M$ and a smooth function $f$ on $M$ satisfying the following properties.
\begin{enumerate}
\item The Reeb graph $W_f$ of $f$ is isomorphic to $G${\rm :} we can take an isomorphism $\phi:W_f \rightarrow G$.
\item If we consider the natural quotient map $q_f:M \rightarrow W_f$ and for each point $\phi(p) \in G$ {\rm (}$p \in W_f${\rm )} that is not a vertex and that is in an edge an integer $q \geq 0$ is assigned to, then the preimage ${q_f}^{-1}(p)$ is a closed and orientable surface of genus $q \geq 0$.
\item For a point $p \in M$ mapped by the quotient map to a vertex $q_f(p) \in W_f$, $f(p)=g \circ \phi(q_f(p))$.
\item At each singular point, the function is represented as a Morse function, a composition of a submersion with a Morse function, or a composition of a fold map with a Morse function. Moreover, for the latter two cases, the Morse functions are height functions of a line or an open ball of dimension $2$. 
\item Furthermore, we can construct $f$ so that the restriction to the preimage of the interior of a small regular neighborhood of vertex of degree $1$ of $W_f$ {\rm (}for the map $q_f${\rm )} is represented as a composition of a fold map with a height function of an open ball of dimension $2$. Furthermore, the fold map is constructed as a map {\rm (}as presented in FIGURE \ref{fig:5}{\rm)} satisfying the following three.
\begin{enumerate}
\item Each connected component of the preimage of each point in the target is a circle, a bouquet of two circles, or a $1$-dimensional polyhedron obtained in the following way.
\begin{enumerate}
\item Choose a suitable integer $l>2$ and take $l$ copies of a circle.
\item For two of the copies, choose exactly one point for each copy, denoted by $p_1$ and $p_{2l-2}$. 
\item For the remaining $l-2$ copies, choose exactly $2$ distinct points for each copy and for the $j$-th copy, denoted by $p_{2j}$ and $p_{2j+1}$.
\item Identify $p_{2j-1}$ and $p_{2j}$ for each $1 \leq j \leq l-1$. 
\end{enumerate}

\item For the singular set of the fold map, remove finitely many singular points and consider the restriction there. This is an embedding.
\item The finitely many singular points removed before are in a same preimage.  
\end{enumerate}
\end{enumerate}
\end{Thm} 

We prove this in the next section. We construct a desired function by the presented method. A new ingredient in the proof is construction around vertices at which the function has local extrema, or STEP 2.
We end the paper by presenting related new problems or Problems \ref{prob:2} and \ref{prob:3}, Remark \ref{rem:3}, explaining relation between studies of the present paper and those in \cite{saeki4}, and Remark \ref{rem:4}.

\section{Several terminologies and the proof of Theorem \ref{thm:1}}
Throughout the paper, ${\mathbb{R}}^k$ denotes the $k$-dimensional Euclidean space. For $x \in {\mathbb{R}}^k$, $||x||$ denotes the distance between the origin $0\in {\mathbb{R}}^k$ and $x$ where the
 underlying metric is the Euclidean metric. $S^k:=\{x \in {\mathbb{R}}^{k+1} \mid ||x||=1.\}$ denotes the $k$-dimensional unit sphere and $D^k:=\{x \in {\mathbb{R}}^{k} \mid ||x|| \leq 1.\}$ denotes the $k$-dimensional unit disc.

For a smooth map, the {\it singular value} is a point such that the preimage contains a singular point. We define the set of all singular values of the map as the {\it singular value set} and the {\it regular value set} as the complementary set of the singular value set. A point in the regular value set of the map is said to be a {\it regular value}, which is explained before. {\it Morse} functions are well-known but we will refer to them in the presentation of a {\it fold} map later. We also note that {\it height functions} of (open) disks are simplest examples of Morse functions and we consider such functions in various cases in this paper.

A smooth map from a manifold of dimension $m>0$ with no boundary into a manifold of dimension $n>0$ with no boundary is said to be a {\it fold} map if the relation $m \geq n$ holds and at each singular point $p$, there exists an integer satisfying $0 \leq i(p) \leq \frac{m-n+1}{2}$ and the map is represented as $(x_1,\cdots,x_m) \mapsto (x_1,\cdots,x_{n-1},\sum_{k=n}^{m-i(p)}{x_k}^2-\sum_{k=m-i(p)+1}^{m}{x_k}^2)$ for suitable coordinates. The case where the target is $\mathbb{R}$ is for {\it Morse} functions. So-called natural {\it height functions} of (open) disks before are Morse functions with exactly one singular point $p$ with $i(p)=0$. Note also that a correspondence of handle attachments to level sets of Morse functions and singular points of Morse functions is well-known. This fact gives us keys in various situations in geometry including several arguments of the present paper. Related to this and fundamental arguments on Morse functions, see \cite{milnor} for example. 
 
\begin{Prop}
For a fold map just before, the following two properties hold.
\begin{enumerate}
\item The integer $i(p)$ is unique {\rm (}we call $i(p)$ the {\it index} of $p${\rm )}. 
\item The set of all singular points of an index is a closed smooth submanifold of dimension $n-1$ with no boundary and the map obtained by restricting the original map to the set is an immersion.  
\end{enumerate}
\end{Prop}

\begin{figure}
\includegraphics[width=40mm]{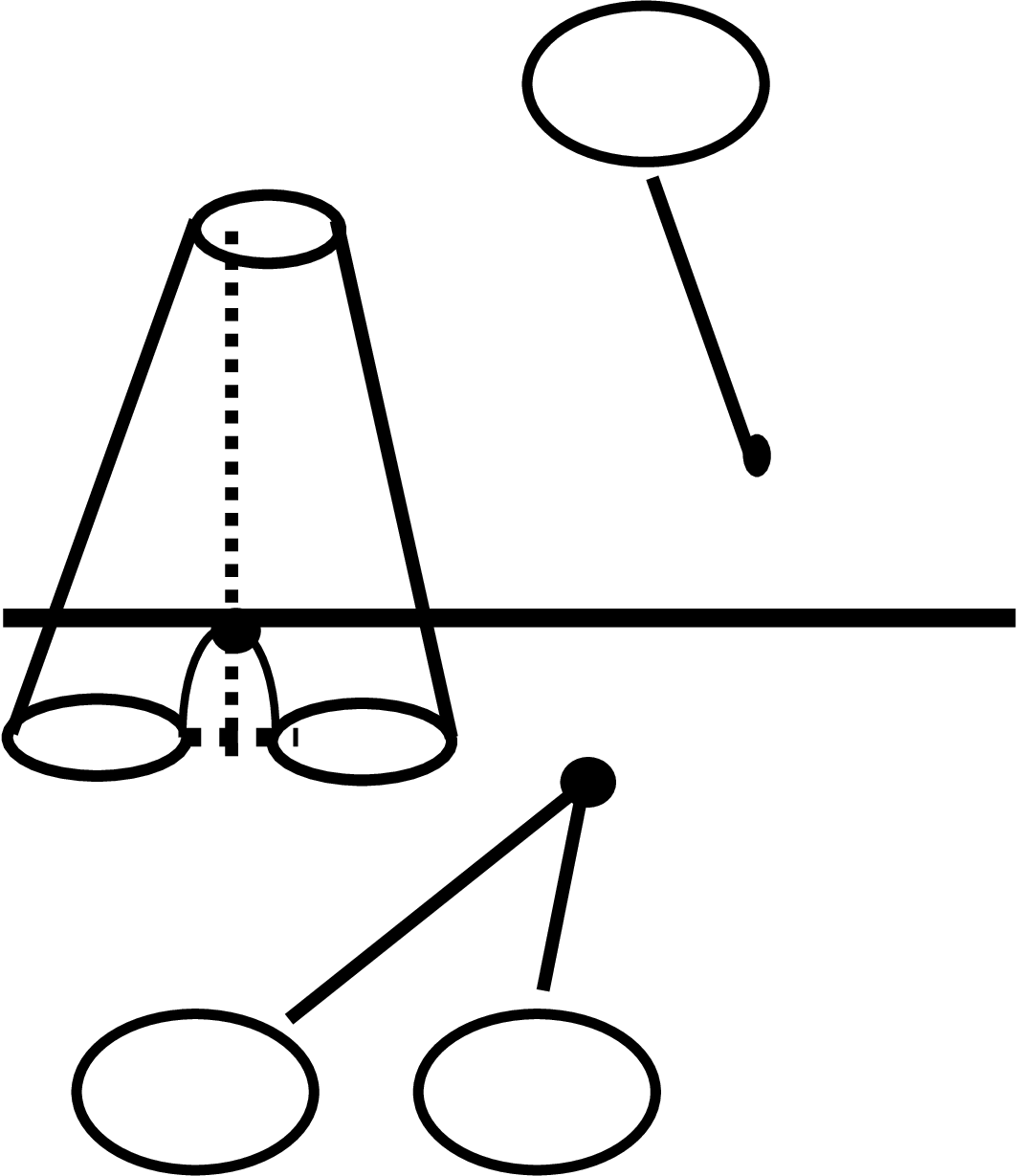}
\caption{The image and the preimage of several points in the target of a fold map on a $3$-dimensional closed manifold into the plane and around its singular values. The horizontal thick line indicates the values at singular points of index $1$ and the given fold map is locally regarded
 as a product of a Morse function on (the interior of) a compact orientable surface of genus $0$ whose boundary (end) consists of three connected components and the identity map on an open interval. Circles indicate the preimages of the regular values pointed by the straight lines.}
\label{fig:1}
\end{figure}

In other words, a fold map is locally regarded as a product of a Morse function with exactly one singular point and the identity map on an open ball. We can also see this by virtue of a kind of Thom's second isotopy lemma (see the proof of Theorem 5.1 of \cite{saekitakase} and so on for explicit usages of such a lemma to see that a given smooth map is regarded as a product map of another smooth map and a suitable identity map on an open disk). See also FIGURE \ref{fig:1}.

For fundamental singular theoretical and differential topological properties of fold maps, see \cite{golubitskyguillemin}, \cite{saeki}, \cite{saeki2} and so on.

We introduce the rigorous definition of the {\it Reeb space {\rm (}graph{\rm )}} of a smooth function. 
Let $X$ be a smooth manifold of dimension $m>1$ and $c:X \rightarrow \mathbb{R}$ be a smooth function.
We define a relation on $X$ by the following rule: $x_1 {\sim}_c x_2$ if and only if $x_1 \in X$ and $x_2 \in X$ are in a same connected component of a preimage $c^{-1}(y)$. This is an equivalence relation on $X$ and we can define the quotient space $W_c:=X/{\sim}_c$. 

\begin{Def}
$W_c$ is said to be the {\it Reeb space} of $c$.
\end{Def} 
We denote the quotient map by $q_c:X \rightarrow W_c$. 
 
We define $V_c:=\{p \mid {q_c}^{-1}(p)$ contains some singular point of $c$. $\}$. 

\begin{Def}
If we can regard $W_c$ as a graph whose vertex set is $V_c$, then $W_c$ is said to be the {\it Reeb graph} of $c$.
\end{Def} 

\begin{Rem}
\label{rem:1}
Reeb spaces may not be regarded as graphs generally (\cite{saeki4}): the present paper concentrates on smooth functions with finitely many singular values and Reeb spaces are Reeb graphs.
\end{Rem}
\begin{Rem}
\label{rem:2}
Most of functions obtained in the present paper are of the class of {\it Morse-Bott} functions (if we remove finitely many singular points). For Morse-Bott functions see \cite{bott} for example.
\end{Rem}

We prove 
Theorem \ref{thm:1}. Hereafter, for a finite set $X$, $\sharp X$ denotes the cardinality of $X$.
\begin{proof}[A proof of Theorem \ref{thm:1}]

\noindent STEP 1 Construction around a vertex $p$ at which $g$ does not have a local extremum. \\
\indent We demonstrate construction presented in \cite{kitazawa2}, which is also based on methods in \cite{masumotosaeki}, \cite{michalak} and \cite{sharko}.
We will construct a local Morse function with exactly one singular value respecting $g$ so that preimages of points in the regular value set are desired connected and orientable surfaces and that the Reeb space is as desired. We introduce notation.

\begin{enumerate}
\item Let $E_1 \neq \emptyset$ be the set of all edges containing the vertex $p$ such that on each edge $g$ has the value smaller than or equal to $g(p)$. Let $S$ be the disjoint union of $\sharp E_1$ closed, connected and orientable surfaces such that the genus of each surface and the integer attached to the corresponding edge coincide.    
\item Let $E_2 \neq \emptyset$ be the set of all edges containing the vertex $p$ such that on each edge $g$ has the values greater than or equal to $g(p)$. Let $S^{\prime}$ be the disjoint union of $\sharp E_2$ closed, connected and orientable surfaces such that the genus of each surface and the integer attached to the corresponding edge coincide.   
\end{enumerate}

This covers all cases we need.

First we consider the case where the numbers assigned to the edges are all $0$. For real numbers $s_1$ and $s_2$ satisfying $s_1<g(p)<s_2$, we can construct a desired function
${\tilde{f}}_{E_1,E_2,p}$ and a suitable homeomorphism ${\phi}_p$ from the resulting Reeb space $W_{{\tilde{f}}_{E_1,E_2,p}}$ of the local function ${\tilde{f}}_{E_1,E_2,p}$ to a small regular neighborhood of $p \in G$ satisfying the following five by virtue of section 4 (Lemma 4.1) of \cite{michalak} together with fundamental theory of Morse functions and handle attachments.
\begin{enumerate}
\item The image of the function is $[s_1,s_2]$ and this function has finitely many or $k>0$ singular points and singular values are always $g(p)$.
\item On the interior they are Morse functions with finitely many singular points.
\item $W_{{\tilde{f}}_{E_1,E_2,p}}$ and the small regular neighborhood of $p \in G$ are regarded as graphs having exactly $\sharp E_1+\sharp E_2$ edges, which are all originating from $p$ and ${\phi}_p(p)$ respectively, and exactly $\sharp E_1+\sharp E_2+1$ vertices including exactly $\sharp E_1+\sharp E_2$ vertices of degree $1$. 
\item The preimage of $s_1$ is diffeomorphic to $S$ and that of $s_2$ is diffeomorphic to $S^{\prime}$.
\item The preimage of the singular value $g(p)$ is homeomorphic to a connected polyhedron obtained in the following way. Take finitely many copies of $S^2$ and let the number of copies. Take a suitable integer $l_0>0$ and disjoint $2l_0$ points there. We divide the set of these $2l_0$  points into $l_0$ pairs. We identify the points for each pair of the points.
\end{enumerate}

Let $s_0$ be a real number satisfying $s_1<s_0<s_2$. Before we consider general cases, we need to define smooth functions $t_{i,s_1,s_0,s_2}$ and $t_{d,s_1,s_0,s_2}$ on a $3$-dimensional compact, connected and orientable manifold whose boundary is diffeomorphic to the disjoint union of a copy of the $2$-dimensional unit sphere $S^2$ and a copy of $S^1 \times S^1$ satisfying the following properties. These functions are compositions of a copy of a Morse function with suitable diffeomorphisms on the targets and the characters ''$i$'' and ''$d$'' in the functions mean ''increasing'' and ''decreasing'' respectively: the orientations for the functions are mutually upside down.

\begin{enumerate}
\item The images of these functions are $[s_1,s_2]$ and these functions have exactly one singular point.
\item On the interior they are Morse functions with exactly one singular point and the singular values are both $s_0$.
\item ${t_{i,s_1,s_0,s_2}}^{-1}(s_1)$ and the connected component of the boundary diffeomorphic to $S^2$ agree. ${t_{d,s_1,s_0,s_2}}^{-1}(s_2)$ and the connected component of the boundary diffeomorphic to $S^2$ agree. 
\item ${t_{i,s_1,s_0,s_2}}^{-1}(s_2)$ and the remaining connected component of the boundary, diffeomorphic to $S^1 \times S^1$, agree. ${t_{d,s_1,s_0,s_2}}^{-1}(s_1)$ and the remaining connected component of the boundary, diffeomorphic to the torus, agree.
\end{enumerate}

We can remove the interiors of two disjoint compact and smooth manifolds the restrictions of the original functions to which give trivial smooth bundles over $[s_1,s_2]$ whose fibers are diffeomorphic to $D^2$ by virtue of (a relative version of) Ehresmann's fibration theorem (\cite{ehresmann}). After that, we can have functions ${\tilde{t}}_{i,s_1,s_0,s_2}$ and ${\tilde{t}}_{d,s_1,s_0,s_2}$ satisfying the following four properties.    
\begin{enumerate}
\item They are restrictions of the original functions.
\item Their images are both $[s_1,s_2]$ and their Reeb spaces are both homeomorphic to a closed interval.
\item They have exactly one singular point and the singular values are both $s_0$.
\item ${{\tilde{t}}_{i,s_1,s_0,s_2}}^{-1}(s_1)$ and ${{\tilde{t}}_{d,s_1,s_0,s_2}}^{-1}(s_2)$ are diffeomorphic to $D^2$.
\item ${{\tilde{t}}_{i,s_1,s_0,s_2}}^{-1}(s_2)$ and ${{\tilde{t}}_{d,s_1,s_0,s_2}}^{-1}(s_1)$ are diffeomorphic to a compact, connected and orientable surface of genus $1$ with one hole (or the torus with one hole).
\end{enumerate}
See also FIGURE \ref{fig:2}. 
\begin{figure}
\includegraphics[height=70mm,width=70mm]{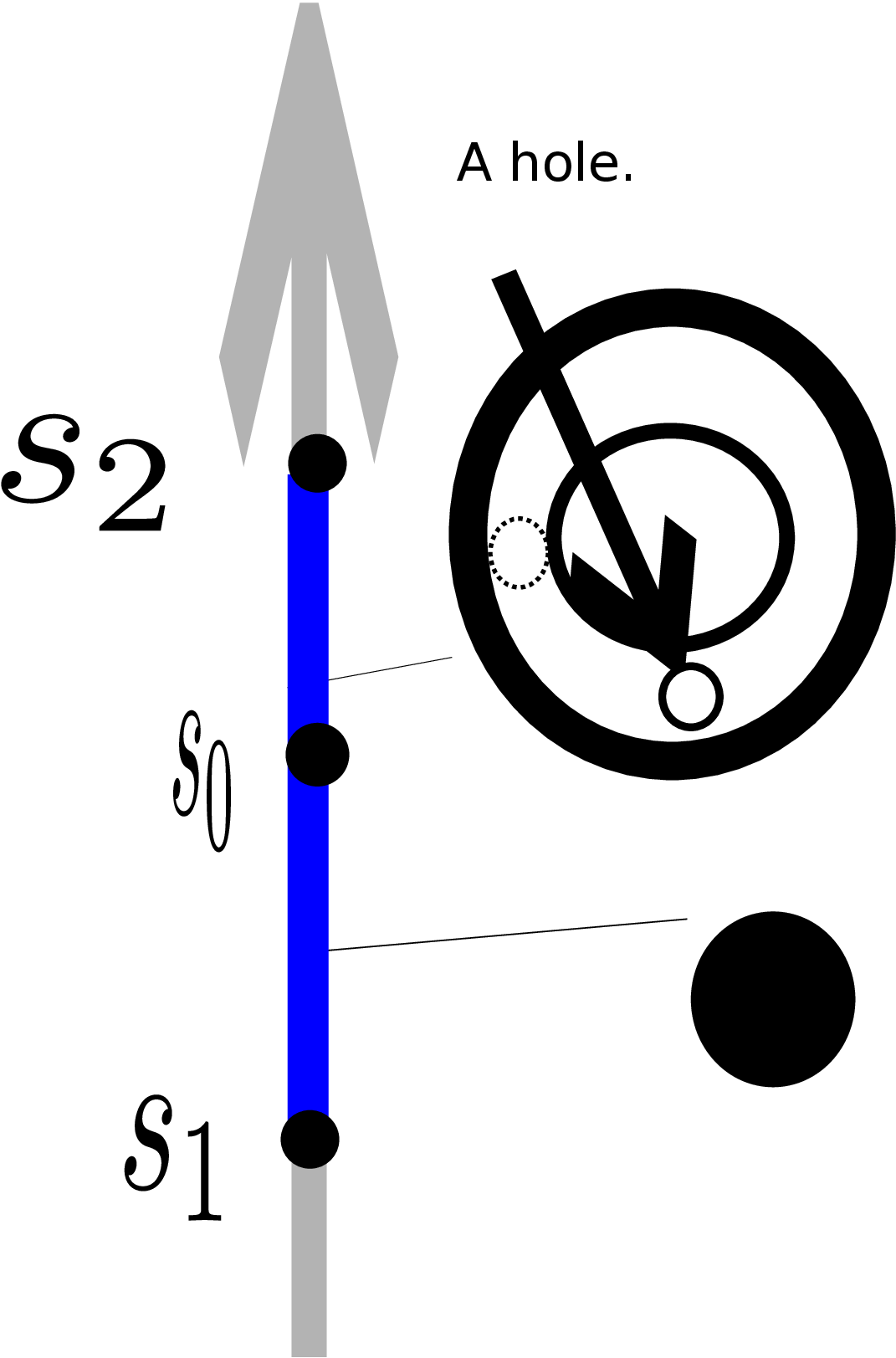}
\caption{The image $[s_1,s_2]$ of the function ${\tilde{t}}_{i,s_1,s_0,s_2}$ (in blue). Preimages (the disk and the torus with a hole in black).}
\label{fig:2}
\end{figure}
We will discuss a general case. For an arbitrary pair $(l_1,l_2)$ of positive integers, we can construct the desired function in the case $(\sharp E_1,\sharp E_2)=(l_1,l_2)$.
 
We construct a desired function ${\tilde{f}}_{E_1,E_2,p}$ by removing the interiors of finitely many disjoint compact and smooth manifolds the restrictions of the function to which give trivial smooth bundles over $[s_1,s_2]=[g(p)-\epsilon,g(p)+\epsilon]$ whose fibers are diffeomorphic to $D^2$ and attach finitely many copies of ${\tilde{t}}_{i,g(p)-\epsilon,g(p),g(p)+\epsilon}$ or ${\tilde{t}}_{d,g(p)-\epsilon,g(p),g(p)+\epsilon}$ for a small $\epsilon>0$ preserving the values of the functions instead.

We can do this and take a homeomorphism ${\phi}_p$ from the resulting Reeb space $W_{{\tilde{f}}_{E_1,E_2,p}}$ of the local function ${\tilde{f}}_{E_1,E_2,p}$ to a small regular neighborhood of $p \in G$ so that the following four properties hold: essentially (a relative version of) Ehresmann's fibration theorem and the structures of the original $3$-dimensional manifold and the original function enable us to do this. Let the function on $E_i$ corresponding the integers assigned to the edges denote by $e_{E_i}$.
\begin{enumerate}
\item $W_{{\tilde{f}}_{E_1,E_2,p}}$ and the regular neighborhood of $p \in G$ are regarded as graphs having exactly $\sharp E_1+\sharp E_2$ edges, which are all originating from $p$ and ${\phi}_p(p)$ respectively, and exactly $\sharp E_1+\sharp E_2+1$ vertices including exactly $\sharp E_1+\sharp E_2$ vertices of degree $1$.
\item The number of copies of ${\tilde{t}}_{i,g(p)-\epsilon,g(p),g(p)+\epsilon}$ we use to construct the desired function is ${\Sigma}_{e \in E_2} e_{E_2}(e)$.
\item The number of copies of ${\tilde{t}}_{d,g(p)-\epsilon,g(p),g(p)+\epsilon}$ we use to construct the desired function is ${\Sigma}_{e \in E_1} e_{E_1}(e)$.
\item The number of the copies of the subset ${{\tilde{t}}_{i,g(p)-\epsilon,g(p),g(p)+\epsilon}}^{-1}([g(p),g(p)+\epsilon])$ of the domain of ${\tilde{t}}_{i,g(p)-\epsilon,g(p),g(p)+\epsilon}$ above mapped to an edge $e^{\prime}$ by $q_{{\tilde{f}}_{E_1,E_2,p}}$ such that ${\phi}_p(e^{\prime})$ is an edge of the small regular neighborhood of $p \in G$ and also in a unique edge $e \in E_2$ of the graph $G$ is $e_{E_2}(e)$.
\item The number of the copies of the subset ${{\tilde{t}}_{d,g(p)-\epsilon,g(p),g(p)+\epsilon}}^{-1}([g(p)-\epsilon,g(p)])$ of the domain of ${\tilde{t}}_{d,g(p)-\epsilon,g(p),g(p)+\epsilon}$ above mapped to an edge $e^{\prime}$ by $q_{{\tilde{f}}_{E_1,E_2,p}}$ such that ${\phi}_p(e^{\prime})$ is an edge of the small regular neighborhood of $p \in G$ and also in a unique edge $e \in E_1$ of the graph $G$ is $e_{E_1}(e)$.
\end{enumerate}

This completes the exposition of a desired function ${\tilde{f}}_{E_1,E_2,p}$ around $p$.

Of course by the construction the five properties in Theorem \ref{thm:1} hold. Note that we may ignore the last property in STEP 1.

This completes STEP 1. \\
\ \\
Hereafter, in the local construction as before we omit  precise expositions as before using a suitable homeomorphism ${\phi}_p$ from the resulting Reeb space of the local function to a small regular neighborhood of $p \in G$. \\
\ \\
\noindent STEP 2 Construction around a vertex $p$ at which $g$ has a local extremum. \\
CASE 1 The case where the vertex $p$ is of degree larger than $1$. \\
\indent We divide the set of all the edges containing the vertex into two non-empty sets $E_1$ and $E_2$. Let $S_i$ be a closed surface obtained as the disjoint union of all closed, connected and orientable surfaces the genius of each of which is the number assigned to the corresponding edge in $E_i$. For $S_1$ and $S_2$, we construct a suitable local Morse function as STEP 1 so that preimages of regular values smaller than the singular value are diffeomorphic to $S_1$ and those of regular values greater than the value are diffeomorphic to $S_2$. After this construction, we embed the image so that the image is a parabola in the plane ${\mathbb{R}}^2$ and that the singular value is the local extremum $g(p)$ with respect to the vertical component and compose a canonical projection to the vertical axis. The resulting local function is a desired function. We can see that this local function satisfies the first three properties in Theorem \ref{thm:1}. For the fourth property, it is regarded as a composition of two Morse functions. At each singular point, the local function is represented as a composition of a submersion with a Morse function or a composition of a Morse function with a Morse function. In addition, in the two compositions, the latter Morse functions are natural height functions of a line. 

For the construction, see also FIGURE \ref{fig:3}.

\begin{figure}
\includegraphics[width=55mm]{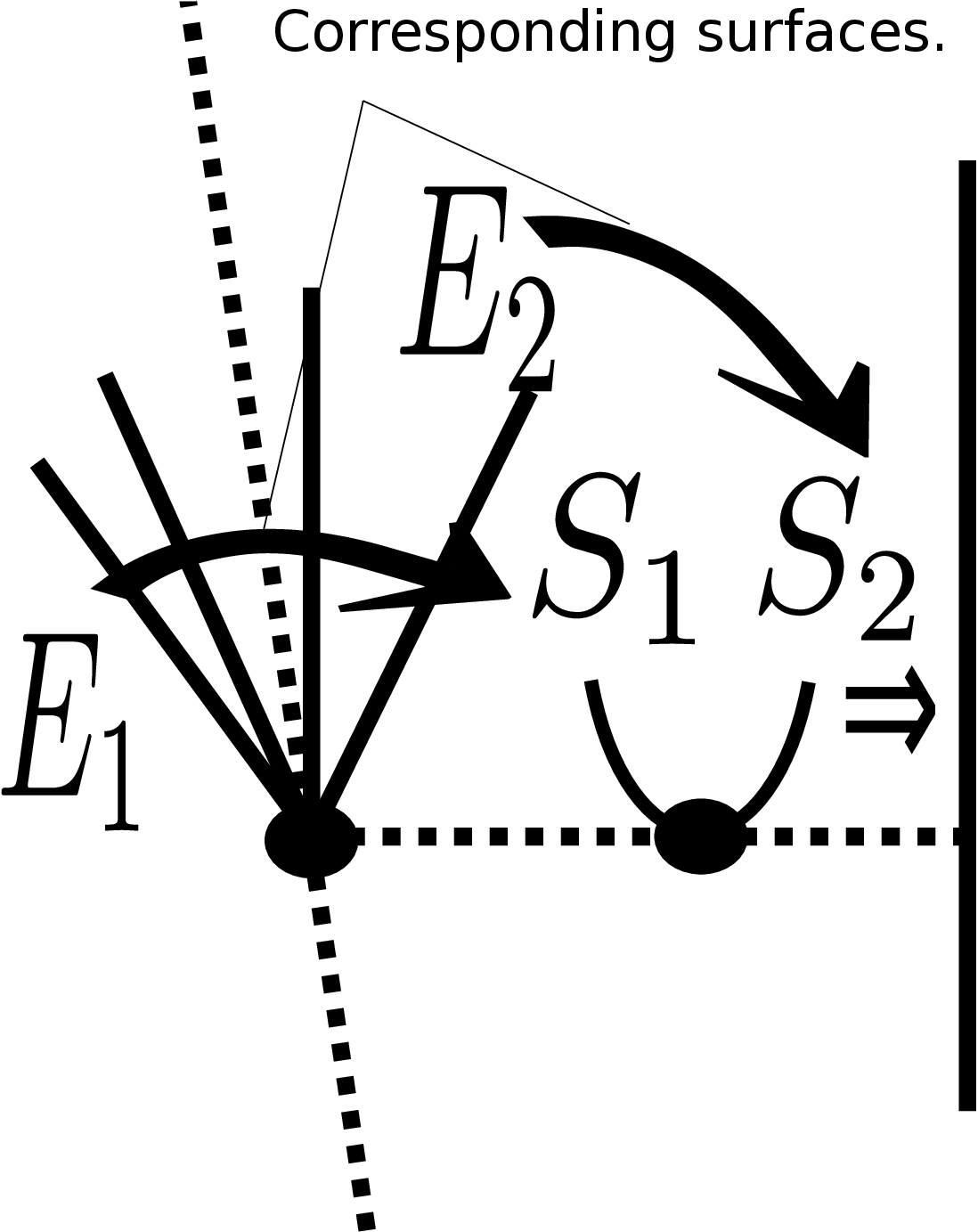}
\caption{Local construction of a function in CASE 1 from a given graph. The left FIGURE shows the vertex and edges containing the vertex, the middle one shows the image of the composition of a local Morse function with an embedding into the plane (and the image is a parabola as explained) and the arrow in the right represents the projection onto the vertical axis. $S_1$ and $S_2$ in the middle FIGURE represent surfaces obtained from the sets $E_1$ and $E_2$ and the assigned numbers and preimages of the (corresponding) points in the image.}
\label{fig:3}
\end{figure}

\ \\
\noindent CASE 2 The case where the vertex $p$ is of degree $1$. \\
\indent We first show the case where a positive integer is assigned to the edge containing the vertex.

Let $k \geq 0$ be an integer, let $a>0$ be a real number and let $\{s_j\}_{j=1}^k$ be an increasing sequence of real numbers satisfying $s_1s_k \leq 0$: if $k=1$ then $s_1=0$ and if $k>1$ then $s_1<0$ or $s_k>0$. 
Let $P$ be a compact, connected and orientable surface obtained by removing the interiors of disjointly and smoothly embedded $k+2$ copies of a $2$-dimensional standard closed disk in a $2$-dimensional sphere.
We have a smooth deformation of smooth functions $F:P \times [-a,a] \rightarrow [-a,a] \times [-a,a]$ defined as $F(b,t):=(f_t(b),t)$ where $\{f_t(b):=F(b,t)\}_{t \in [-a,a]}$ are smooth functions on $P$ satisfying the following properties. 
\begin{enumerate}
\item The images of these functions are $[-a,a]$ and these functions have exactly $k$ singular points.
\item On the interior they are Morse functions with exactly $k$ singular points whose indices are $1$ and the singular value set of $f_t$ is $\{s_j-\frac{s_j}{a}(t+a) \mid 1 \leq j \leq k\}$.
\item For $f_t$, the preimage of a regular value $p \in [-a,a]$ in the lower bound of its singular value set is a circle and the preimage of a regular value $p \in [-a,a]$ in the upper bound of the singular value set is a disjoint union of $k+1$ circles.
\item For $-a<t<0$ and $f_t$, the preimage of a regular value in $\{p \in [-a,a] \mid s_j-\frac{s_j}{a}(t+a)<p<s_{j+1}-\frac{s_{j+1}}{a}(t+a)\}$ is a disjoint union of $1+j$ circles. For $0<t<a$ and $f_t$, the preimage of a regular value in $\{p \in [-a,a] \mid s_{j+1}-\frac{s_{j+1}}{a}(t+a)<p<s_{j}-\frac{s_{j}}{a}(t+a)\}$ is a disjoint union of $k+1-j$ circles.
\item The preimage of each singular value of these functions is as follows. 
\begin{enumerate}
\item A disjoint union of finitely many copies of a circle and a bouquet of two circles in the case $f_t$ ($t \neq 0$).
\item A $1$-dimensional polyhedron obtained in the following way in the case $f_t$ ($t=0$).
\begin{enumerate}
\item Take $k+1$ circles.
\item For two of the copies, choose exactly one point for each copy, denoted by $p_1$ and $p_{2k}$. 
\item For the remaining $k-1$ copies, choose exactly $2$ distinct points for each copy and for the $j$-th copy, denoted by $p_{2j}$ and $p_{2j+1}$.
\item Identify $p_{2j-1}$ and $p_{2j}$ for each $1 \leq j \leq k$. 
\end{enumerate}
\end{enumerate}
\item (The restriction of) $F$ (to the interior) is a fold map and the singular set of $F$ is a disjoint union of $k$ copies of a closed interval. If we restrict $F$ to each connected component of the singular set, then the resulting map is an embedding and the image is $\{(s_j-\frac{s_j}{a}(t+a),t) \in [-a,a] \times [-a,a] \mid -a \leq t \leq a, 1 \leq j \leq k\}$.
\end{enumerate}   

See also FIGURE \ref{fig:4}. Note that this makes sense for $k=0$.

\begin{figure}
\includegraphics[width=35mm]{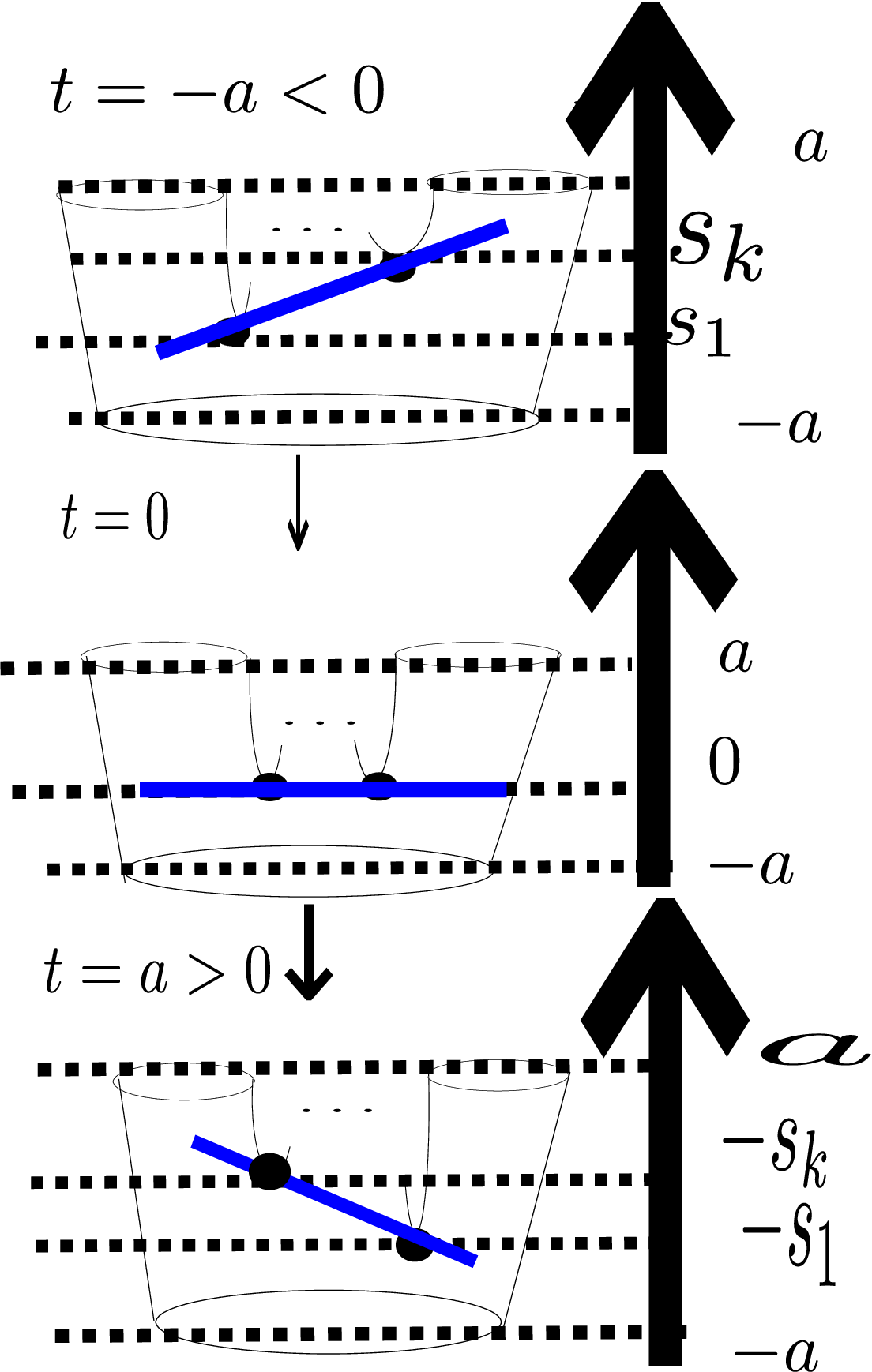}
\caption{Smooth functions $f_t(b)$ on $P$. This also represents a deformation $F$ of such functions on $t \in [-a,a]$. Dots are singular points.}
\label{fig:4}
\end{figure}

Thus we have a local fold map from a $3$-dimensional orientable manifold onto $(-a,a) \times (-a,a) \subset {\mathbb{R}}^2$ (by restricting the domain and the target of $F$ to the interiors). FIGURE \ref{fig:5} shows the image of the fold map.
Thick lines intersect at the origin $0 \in {\mathbb{R}}^2$ and these lines represent straight lines containing the origin $0 \in {\mathbb{R}}^2$ consisting of values of the map $F$ at singular points whose indices are $1$. If we go in the direction of the arrows in the plane, then the number of connected components of the preimages of regular values increases. 

We restrict the map to the preimage of the disk bounded by the dotted circle in FIGURE \ref{fig:5} and let the target be this disk. 
We take the disk as $\{x \in {\mathbb{R}}^2 \mid ||x|| \leq r\}$ for a (sufficiently small) positive number $r>0$.
Let the map on the $3$-dimensional manifold ${\bar{M}}_p$ onto the disk denote by $\tilde{F_p}$.


The restriction of $\tilde{F_p}$ to the preimage of $\{x \in {\mathbb{R}}^2 \mid ||x||= r\}$ is a so-called {\it circle-valued Morse} map and this has exactly $2k$ singular points. Their indices are all $1$. We can see that the preimage of the embedded dotted circle $\{x \in {\mathbb{R}}^2 \mid ||x||= r\}$ is a closed, connected and orientable surface of genus $k+1$. Note also that in cases where the numbers of the thick lines are $0$ or $1$ (or $k=0$ or $k=1$) or so-called generic cases, see \cite{saeki2} for example. The local map in a general case is also regarded as a generalization of these maps.

\begin{figure}
\includegraphics[width=35mm]{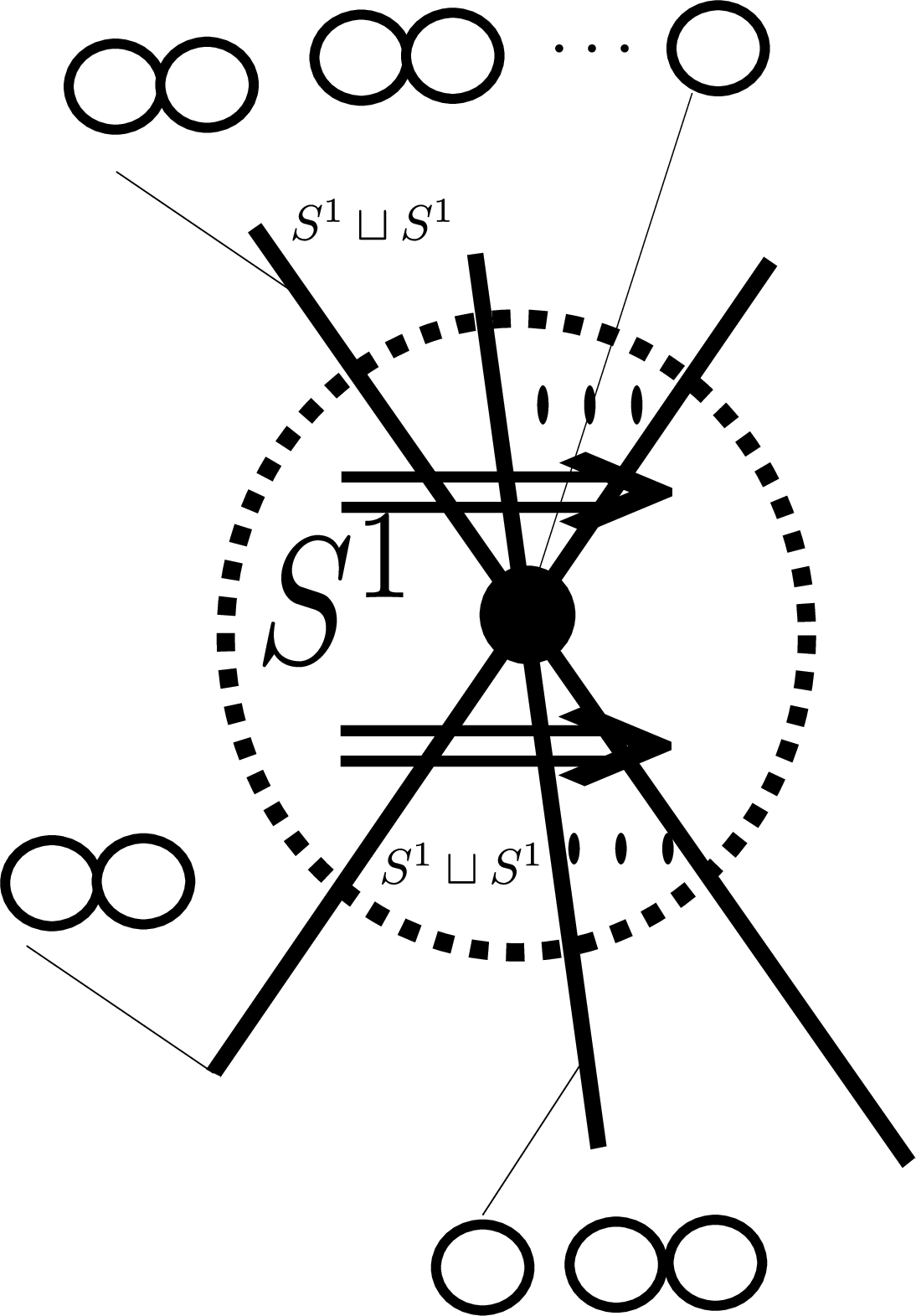}
\caption{The image (the singular value set) of a local fold map into the plane: the manifolds such as $S^1$ and $S^1 \sqcup S^1$ and the shapes consisting of circles represent preimages of the corresponding points.}
\label{fig:5}
\end{figure}
We consider a natural height function of the (interior of the) disk $\{x \in {\mathbb{R}}^2 \mid ||x|| \leq r\}$ whose exactly one singular point is the origin $0$ (we can define the height function as a function represented
as $(x_1,x_2) \mapsto \pm({x_1}^2+{x_2}^2)+g(p)$ for suitable coordinates). 

By composing the local fold map $\tilde{F_p}$ ($\tilde{F_p} {\mid}_{{\tilde{F_p}}^{-1}(\{x \in {\mathbb{R}}^2 \mid ||x|| < r\})}$) on the $3$-dimensional manifold into the (interior of the) disk $\{x \in {\mathbb{R}}^2 \mid ||x|| \leq r\}$ with this height function, we obtain a local function such that the image and the Reeb space is homeomorphic to $[0,1)$. Each point in the interior of the image is a regular value since the image of the rank of the differential of the fold map at each singular point is $1$ and the preimage is a closed, connected and orientable surface of genus $k+1$. For this, note also that for the fold map into the plane (disk), the preimage of $\{x \in {\mathbb{R}}^2 \mid ||x||=t>0\}$ is a closed, connected and orientable surface of genus $k+1$ and that the image of the differential at each singular point mapped to a point or a singular value different from the origin $0 \in {\mathbb{R}}^2$ is regarded as the unique straight line in the singular value set of the fold map containing the singular value. The resulting local function is, at each point in the preimage of $0 \in [0,1)$, represented as a composition of a fold map into the plane with a height function of an open $2$-dimensional disk. 

This local function is a desired local function. For the functions and maps here, we can see that the fourth and fifth properties in Theorem \ref{thm:1} hold. 
For the fold map into the plane, the preimage of the origin $0 \in {\mathbb{R}}^2$ may contain more than one singular point, where other preimages do not. 
As in CASE 1, we can see that this local function satisfies the first three properties in the theorem.

We must also show the case where $0$ is assigned to the edge containing the vertex. 
To make the preimage at the point in the interior of the edge diffeomorphic to $S^2$, or equivalently, a closed, connected and orientable surface of genus $0$, it is sufficient to consider a height function of an open $2$-dimensional disk.

This completes the proof for CASE 2 and STEP 2.  \\
 \\
\noindent STEP 3 Around each connected component of the complementary set of the disjoint union of small regular neighborhoods of vertices. \\
\indent Last, we construct functions around each connected component of the complementary set of the disjoint union of small regular neighborhoods of vertices. We can construct these functions as trivial smooth bundles. \\
 \ \\
\indent Finally, we can glue the local functions together by suitable diffeomorphisms. Note that several local functions in the present paper are constructed as functions on open manifolds with compact ends and we can extend them to smooth functions on compact manifolds with boundaries in canonical ways. Note also that we can choose the diffeomorphisms so that the manifolds are orientable. This gives a desired function $f$ on a $3$-dimensional closed, connected and orientable manifold and we also have a desired isomorphism $\phi$ between the graphs by gluing ${\phi}_p$ in a natural way. This completes the proof.  
\end{proof}
 
\begin{Prob}
\label{prob:2}
Find restrictions on the topologies of manifolds admitting functions of Theorem \ref{thm:1} satisfying appropriate topological conditions.
\end{Prob}

As a simplest answer, in \cite{saeki3}, manifolds admitting Morse functions such that preimages of regular values are disjoint unions of a copy of $S^k$ have been studied for $k \geq 1$. As another simplest result, it has been shown that a $3$-dimensional closed, connected and orientable manifold admits a Morse function such that connected components of preimages of regular values are always diffeomorphic to $S^2$ or $S^1 \times S^1$ if and only if it is diffeomorphic to $S^3$ or represented as a connected sum of Lens spaces or copies of $S^2 \times S^1$. In other words, the case where $0$ or $1$ are assigned to all edges of the graphs in Theorem \ref{thm:1} was studied there under the additional constraint that the class of functions is restricted to the class of Morse functions. See also Remark \ref{rem:3}.

\begin{Prob}
\label{prob:3}
Can we show a non-orientable version of Theorem \ref{thm:1}?
\end{Prob}

The author has a positive idea and an answer will be presented in a forthcoming paper.

\begin{Rem}
\label{rem:3}
\cite{saeki4} is regarded as a related work motivated by the present paper. This concerns a general case where preimages of regular values are arbitrary compact manifolds possibly with boundaries of a fixed dimension. In the work explicit singularities of smooth functions are not studied: only so-called bump functions are used in constructing global functions for example. In the present paper an explicit class of smooth functions on $3$-dimensional closed, connected and orientable manifolds is studied where such very explicit cases are not studied in \cite{saeki4}. 

Note that for example if we restrict our functions to functions with finitely many singular points then they are Morse functions whose singular values may agree at distinct singular points. Related to this, note also that in the answer to Problem \ref{prob:2} presented just before, such a restriction of the class of functions restrict the class of manifolds admitting functions of the class to the presented class of $3$-dimensional manifolds. Furthermore, if we do not restrict the class of functions of Theorem \ref{thm:1} and we pose the condition that to all edges of the given graph $G$ $0$ or $1$ are assigned only, then we do not know the class of manifolds admitting functions of the class. Understanding this class seems to be difficult.
\end{Rem}

\begin{Rem}
\label{rem:4} 
\cite{kitazawa} and \cite{kitazawa2} also present examples of Morse functions of the class of functions in Theorem \ref{thm:1}.
\end{Rem}
\section{Acknowledgement.}
\label{sec:3}
\thanks{The author would like to thank Professor Irina Gelbukh for invaluable comments given after the submission of an earlier version of the present paper to arXiv.org e-Print archive (https://arxiv.org/). One of the comments let the author think Problem \ref{prob:2} as a fundamental and important problem. She also kindly told the author about \cite{gelbukh} and \cite{gelbukh2}. \\
\indent The author is a member of the project JSPS KAKENHI Grant Number JP17H06128 "Innovative research of geometric topology and singularities of differentiable mappings"
 and this work is supported by this project. The author would like to thank Osamu Saeki and colleagues in this project. The author can present several stuffs in the present paper owing to their comments in an informal seminar on the present paper. The author would like to thank Osamu Saeki for publishing \cite{saeki4} as a related work and an informal discussion on this.
 
The author would like to thank anonymous referees, who improved the quality of the present paper.

\end{document}